\documentclass[a4paper,11pt]{article}
\title{Antipodal Hadwiger numbers of finite-dimensional Banach spaces}

\author{S.K.Mercourakis and G.Vassiliadis}

\date{}

\usepackage{amsmath}
\usepackage{amsfonts}
\usepackage{amssymb}
\usepackage{amscd}
\usepackage{amsthm}
\usepackage{makeidx}
\usepackage{euscript}
\usepackage{graphicx}

\theoremstyle{plain}
\newtheorem{theo}{Theorem}
\newtheorem*{theor}{Theorem}
\newtheorem{lemm}{Lemma}
\newtheorem{prop}{Proposition}

\theoremstyle{definition}
\newtheorem{defn}{Definition}
\newcommand{\beg}{\begin{proof}}
\newcommand{\fdb}{finite-dimensional Banach $\,$}
\newcommand{\bsa}{bounded and separated antipodal $\,$}
\newcommand{\ban}{Banach space $\,$}

\begin{document}
\maketitle

\begin{abstract}
\footnotesize

Let $X$ be a finite-dimensional Banach space; we introduce and investigate a natural
generalization of the concepts of Hadwiger number $H(X)$ and strict Hadwiger number
$H'(X)$. More precisely, we define the antipodal Hadwiger number $H_\alpha(X)$ as the largest
cardinality of a subset $S \subseteq S_X$, such that $\forall x \neq y \in S \,\,\,
\exists f \in B_{X^*}$ with 
\[1 \le f(x)-f(y) \,\,\, \textrm{and} \,\,\, f(y) \le f(z) \le f(x) \,\,\, \textrm{for}
\,\,\, z \in S.\]
The strict antipodal Hadwiger number $H'_\alpha(X)$ is  defined analogously.
We prove that $H'_\alpha(X)=4$ for every Minkowski plane and estimate (or in 
some cases compute) the numbers $H_\alpha(X)$ and $H'_\alpha(X)$, where $X=\ell_p^n,
1 < p \le +\infty$ and $n \ge 2$. We also show that the number $H'_\alpha(X)$ grows
exponentially in $\dim X$.

\normalsize
\end{abstract}

\footnote{\noindent 2010 \textsl{Mathematics Subject
Classification}: Primary 46B20;Secondary 52C99.\\
\textsl{Key words and phrases}: Hadwiger number, antipodal set, Banach-Mazur distance  .}

\section*{Introduction}
If $X$ is any (real) Banach space, then $B_X$ and $S_X$ denote its closed unit ball
and unit sphere respectively. A subset $S$ of a normed space $X$ is said to be \textit{$\delta$-separated},
if $\|x-y\| \ge \delta$ for $x \neq y \in S$. Specifically $S$ is called \textit{equilateral}, if there is a $\lambda>0$
such that for $x \neq y \in S$ we have $\|x-y\|=\lambda$; we also call $S$ a $\lambda$-equilateral
set. Any equilateral set in an $n$-dimensional space is of cardinality at most $2^n$
and the maximum is attained only when $X=\ell_\infty^n$ (see \cite{S}).

Let $X$ be a \fdb space. The \textit{Hadwiger number} $H(X)$ of $X$ is the largest cardinality
of a set $S \subseteq S_X$ such that $\|x-y\| \ge 1$, for $x \neq y \in S$. Also the \textit{strict
Hadwiger number} $H'(X)$ of $X$ is the largest cardinality
of a set $S \subseteq S_X$ such that $\|x-y\|> 1$, for $x \neq y \in S$. It is clear that $H'(X) \le H(X)$.
There exists an extensive bibliography on the above concepts (see \cite{SWAN}).

A subset $S$ of a normed space $X$ is said
to be \textit{antipodal} if for every $x,y \in S$ with $x \neq y$ there
exists $f \in X^*$ such that $f(x)<f(y)$ and $f(x) \le f(z) \le
f(y) \, \forall z \in S$. That is, for every $x,y \in S$ with $x
\neq y$ there exist closed distinct parallel support hyperplanes
$P(=\{z \in X:f(z)=f(x)\})$ and $Q(=\{z \in X:f(z)=f(y)\})$ with
$x \in P$ and $y \in Q$. Every antipodal subset of an $n$-dimensional real vector space
has cardinality at most $2^n$ by a result of Danzer and Gr\"{u}nbaum (see \cite{DG}, also \cite{S})
and this is attained only when the points of
the antipodal set are the vertices of an $n$-dimensional
parallelotope.

A \textit{bounded and separated antipodal subset} of a normed space $X$ is a subset $S \subseteq B_X$,
for which there is $d>0$ such that $\forall x \neq y \in S$ there is $f \in B_{X^*}$ with
$d \le f(x)-f(y)$ and $f(y) \le f(z) \le f(x)$ for $z \in S$. If X is a finite-dimensional real
vector space then this concept of antipodality coincides with the
classical one (see \cite{DG}, \cite{P}). The generalization of antipodality stated above was defined in \cite{MV}
where the following Theorem (Th. 3 of \cite{MV}) was proved

\begin{theor}
Let $(X,\| \cdot \|)$ be a \ban and $S \subseteq X$ be a bounded
and separated antipodal set with constant $d$. Then we have:
\begin{enumerate}
\item There is an equivalent norm $||| \cdot |||$ on $X$, such that
S is an equilateral set in $(X,||| \cdot |||)$.
\item The Banach-Mazur distance between $(X,\| \cdot \|)$ and $(X,||| \cdot
|||)$ satisfies the inequality $\,d((X,\| \cdot \|),(X,||| \cdot
|||)) \leq \frac{2}{d}$.
\end{enumerate}
\end{theor} 

The \textit{Banach-Mazur distance} between two isomorphic Banach spaces $X$ and $Y$
is $d(X,Y)=\inf \{\|T\| \cdot \|T^{-1}\|, \text{where } T:X \rightarrow Y \:\text{ is an isomorphism} \}$. It is easy to see that
the dual spaces also have the same Banach-Mazur distance $d(X,Y)=d(X^*,Y^*)$. An equivalent
(geometric) definition of the Banach-Mazur distance is (in finite dimensions), for $K,L \subseteq \mathbb{R}^n$ symmetric convex bodies 
$d(K,L)=\inf\{r>0: L \subseteq T(K) \subseteq r \cdot L,\text{where} \: T: \mathbb{R}^n \rightarrow \mathbb{R}^n \:
\text{is a linear transformation} \}$.

Let $X$ be an $n$-dimensional Banach space. An \textit{Auerbach basis} of $X$ is a biorthogonal system 
$\{(e_i,e^*_i):i=1,2,\dots,n\}$ in $X \times X^*$ (i.e. $e^*_i(e_j)=\delta_{ij}, \, i,j=1,2,\dots,n$)
such that $\{e_i:i=1,2,\dots,n\}$ is a basis of $X$ and $\|e_i\|=\|e^*_i\|=1$ for
$i=1,2,\dots,n$. It is well known that any finite-dimensional Banach space admits an Auerbach basis (see \cite{D}).

In the present paper we introduce and study some interesting analogues of the Hadwiger and
the strict Hadwiger number for a finite-dimensional Banach space, which we call 
antipodal Hadwiger ($H_\alpha(X)$) and strict antipodal Hadwiger number ($H'_\alpha(X)$).
An analogue of the equilateral number $e(X)$, denoted by $e_c(X)$ is defined, which roughly
is the largest cardinality of an equilateral subset of $X$ with center.
The main results are the following:

\begin{enumerate}
\item We prove that $H'_\alpha(X)=4$ for every Minkowski plane (Prop.4, Remarks 3(1)).
\item We estimate and in some cases find exact values of the numbers $H_\alpha(X)$ and 
$H'_\alpha(X)$, when $X=\ell_p^n,1<p \le + \infty$ and $n \ge 2$ (Th.1). We also show that for
an $n$-dimensional Banach space $X$ the number $H'_\alpha(X)$ increases exponentially with $n$ (Th.4)
and the number $e_c(X)$ is bounded below by an unbounded function $\varphi(n)$ (Th.3).
\item We compute the numbers $H_\alpha(X)$ and $H'_\alpha(X)$ when $X=\ell_1^3$ or $X$ is the 
Petty space, i.e. $X=(\mathbb{R}^3,\|\cdot\|)$, where $\|(x,y,z)\|=\sqrt{x^2+y^2}+|z|$ (Prop.10).
\end{enumerate}

\textbf{Acknowledgements.}
The authors wish to thank the anonymous referee for valuable suggestions
and remarks, which helped us to complete our initial investigation and greatly
contributed to the final form of this paper.

\section*{Antipodal Hadwiger number of a \fdb space}

We will always assume that $X$ is a \fdb space. Given the 
definition of Hadwiger number and its variants (see \cite{SWAN}) it is natural
to introduce the following definitions:

\begin{defn}

\begin{enumerate}
\item The \textit{antipodal Hadwiger number} $H_\alpha(X)$ is the largest cardinality
of a set $S \subseteq S_X$ such that $\forall x \neq y \in S$ there is $f \in B_{X^*}$ with
\[ 1 \le f(x)-f(y) \,\, \textrm{and}  \,\, f(y) \le f(z) \le f(x) \,\, \textrm{for} \,\, z \in S. \]
Clearly $S$ is a \bsa set of unit vectors with $d=1$. In particular $S \subseteq S_X$
and $\|x-y\| \ge 1$ for $ x \neq y \in S$, hence $H_\alpha(X) \le H(X)$.

\item The \textit{strict antipodal Hadwiger number} $H'_\alpha(X)$ is the largest cardinality
of a set $S \subseteq S_X$ such that $\forall x \neq y \in S$ there is $f \in B_{X^*}$ with
\[1 < f(x)-f(y) \,\, \textrm{and} \,\, f(y) \le f(z) \le f(x) \,\, \textrm{for} \,\, z \in S. \]
As above, $S \subseteq S_X$
and $\|x-y\| > 1$ for $ x \neq y \in S$, hence $H'_\alpha(X) \le H'(X)$.

\end{enumerate}

\end{defn}

\noindent \textbf{Remarks 1}\\
\noindent (1) Since every antipodal subset of an $n$-dimensional real vector space
has cardinality at most $2^n$, we get that $H'_\alpha(X) \le H_\alpha(X) \le 2^n$.\\
\noindent (2) If $Y$ is a subspace of $X$, then obviously $H_\alpha(Y) \le H_\alpha(X)$
and $H'_\alpha(Y) \le H'_\alpha(X)$.\\
\noindent (3) Let $S$ be an antipodal subset of $S_X$ such that $\|x-y\| \ge 1$
for $ x \neq y \in S$. Since the space is finite-dimensional, $S$ is a \bsa set,
but the constant $d$ (see Definition 1) may be smaller than 1. 

A simple example is a set of three 
consecutive vertices of a regular hexagon inscribed in the unit circle of $\ell_2^2$.
These form an isosceles and obtuse triangle (with an angle of $120^\circ$)  with 2 equal sides of length 1. 
Any functional
$f \in B_{\ell_2^2}$ separating vertices $x,y$ which are at distance 1 gives an evaluation
$|f(x)-f(y)|<1$.\\
\noindent (4) Any $\lambda$-equilateral set $S \subseteq B_X$ is a \bsa set with $d=\lambda$
(see \cite{MV}, Proposition 2). Since the usual basis $S=\{e_1,e_2,\dots,e_n\}$ of $\ell_p^n, 1<p<\infty$ is a 
$2^{\frac{1}{p}}$-equilateral set and $2^{\frac{1}{p}}>1$, it follows in particular that $H'_\alpha(\ell_p^n) \ge n$.
If $p=1$, then the set $\{\pm e_k:k=1,2,\dots,n\}$ is 2-equilateral, hence $H'_\alpha(\ell_1^n) \ge 2n$.

In the sequel we will obtain lower estimates for the antipodal and the strict antipodal Hadwiger 
numbers of a \fdb space.

\begin{prop}
Let X be an $n$-dimensional Banach space. Also let $\{(e_i,e^*_i):i=1,2,\dots,n\}$
be an Auerbach basis of $X$. Then the set $A=\{\pm e_i, i=1,2,\dots,n\}$ is a \bsa
subset of $X$ with constant $d=1$ and hence $H_\alpha(X) \ge 2n$.
\end{prop}

\beg
We check that $\forall x \neq y \in A$ there is $f \in B_{X^*}$ with
$1 \le f(x)-f(y)$ and $f(y) \le f(z) \le f(x)$ for $z \in A$. We have the following cases:\\
\noindent (1) Let $x=e_i$ and $y=e_j$, $i \neq j$. Set $f=\frac{e^*_i-e^*_j}{2}$, then $\|f\| \le 1$ and
\[-\frac{1}{2}=f(e_j) \le f(\pm e_k) \le f(e_i)=\frac{1}{2} \,\, \textrm{for} \,\, k=1,2,\dots,n \]
\noindent (2) Let $x=e_i$ and $y=-e_i$. Set $f=e^*_i$, then $\|f\|=1$ and
\[-1=f(-e_i) \le f(\pm e_k) \le f(e_i)=1 \,\, \textrm{for} \,\, k=1,2,\dots,n \]
\noindent (3) Let $x=e_i$ and $y=-e_j$, $i \neq j$. Set $f=\frac{e^*_i+e^*_j}{2}$, then $\|f\| \le 1$ and
\[-\frac{1}{2}=f(-e_j) \le f(\pm e_k) \le f(e_i)=\frac{1}{2}  \,\, \textrm{for} \,\, k=1,2,\dots,n \]
\noindent (4) Let now $x=-e_i$ and $y=-e_j$, $i \neq j$.
This case is similar to (1) and the proof is complete.
\end{proof}

In case when the space is smooth we have the following stronger result, which was proved in \cite{G} (Prop. 3.9);
see also \cite{GM} (Prop. 1.6):

\begin{prop}
Let X be an $n$-dimensional smooth Banach space and let $\{(e_i,e^*_i):i=1,2,\dots,n\}$
be an Auerbach basis of $X$. Then the set $A=\{\pm e_i, i=1,2,\dots,n\}$ is a \bsa
subset of $B_X$ with constant $d=1+\varepsilon$ for some $\varepsilon>0$. So when $X$ is smooth
we have $H'_\alpha(X) \ge 2n$.
\end{prop}

Strengthening our assumption about the basis (supposing it is 1-suppression unconditional) we can prove
the following:

\begin{prop}
Let X be an $n$-dimensional Banach space and let $\{e_i,1 \le i \le n\}$
be a 1-suppression unconditional normalized basis of $X$. If the norm of $X$ is strictly convex
(or smooth), then the set $A=\{\pm e_i, i=1,2,\dots,n\}$ is a \bsa
subset of $B_X$ with constant $d=1+\varepsilon$, hence $H'_\alpha(X) \ge 2n$.
\end{prop}

\beg
Recall that the basis $\{e_i,1 \le i \le n\}$ of $X$ is 1-suppression unconditional, if for any
$\alpha_1,\alpha_2, \dots,\alpha_n \in \mathbb{R}$ and $F \subseteq \{1,2,\dots,n\}$ we have 
$\|\sum_{k \in F} \alpha_k e_k\| \le \|\sum_{k=1}^n \alpha_k e_k\|$. Since the basis is normalized,
we get that the biorthogonal functionals $\{e^*_i,1 \le i \le n\}$ are also normalized. In
particular $\{(e_i,e^*_i):1 \le i \le n\}$ is an Auerbach basis of $X$.

In both cases we use the fact that $\|e_i \pm e_j\|>1 \,\,\, \text{for } 1 \le i < j \le n$.
In any case we have $\|e_i \pm e_j\| \ge 1$, since $\{(e_i,e^*_i):1 \le i \le n\}$ is an Auerbach basis of $X$.

Let $X$ be strictly convex and assume that $\|e_i+e_j\|=1$ for some $1 \le i < j \le n$. Then 
$e^*_i(e_i+e_j)=e^*_i(e_i)=1$, hence the normalized functional $e^*_i$ attains its maximum at two
distinct points of the unit ball, a contradiction (see \cite{SG} \S 3.2). So $\|e_i+e_j\|>1$ and
similarly $\|e_i-e_j\|>1$.

Let $d=\min\{\|e_i \pm e_j\|: 1 \le i < j \le n \}$. Then $d=1+\varepsilon$ for some $\varepsilon>0$. We
prove the following:\\
\noindent (I) Given $1 \le i < j \le n$ there are $\lambda,\mu \in \mathbb{R}$ such that the functional
$f=\lambda e^*_i+\mu e^*_j$ satisfies $\|f\| \le 1$, $f(e_i-e_j)=\|e_i-e_j\|=\lambda-\mu \ge 1+\varepsilon$
and $\lambda=f(e_i) \ge f(\pm e_k) \ge f(e_j)=\mu$, for $1 \le k \le n$.

To prove this, note that the set $\{e^*_i:1 \le i \le n\}$ is also 1-suppression unconditional normalized 
basis of $X^*$ and thus $\|e_i-e_j\|=\sup\{(xe^*_i+ye^*_j)(e_i-e_j):\|xe^*_i+ye^*_j\| \le 1\}$. So 
there are $\lambda,\mu \in \mathbb{R}$ such that the functional $f=\lambda e^*_i+\mu e^*_j$ gives
\[f(e_i-e_j)=\|e_i-e_j\|=\lambda-\mu \ge d=1+\varepsilon.\]
It follows that \\
\noindent(a) $-1 \le \lambda,\mu \le 1$ (since $f(e_i)=\lambda,f(e_j)=\mu$, $\|f\| \le 1$ and $\|e_i\|=\|e_j\|=1$).\\
\noindent (b) $-1 \le \mu<0<\lambda \le 1$ (since $\lambda-\mu \ge 1+\varepsilon$). \\
\noindent Consequently $\lambda=f(e_i) \ge f(\pm e_k) \ge f(e_j)=\mu$, for $1 \le k \le n$.

\noindent (II) Given $1 \le i < j \le n$ there are $\lambda,\mu \in \mathbb{R}$ such that the functional
$g=\lambda e^*_i+\mu e^*_j$ satisfies $\|g\| \le 1$, $g(e_i+e_j)=\|e_i+e_j\|=\lambda+\mu \ge 1+\varepsilon$
and $-\lambda=g(-e_i) \le g(\pm e_k) \le g(e_j)=\mu$, for $1 \le k \le n$.

The proof of this is similar.

\noindent (III) For $1 \le i \le n$ set $f=e^*_i$. Then we have 
$-1 = f(-e_i) \le f(\pm e_k) \le f(e_i)=1$, for $1 \le k \le n$, which completes the proof that
$A$ is a \bsa subset of $X$.

Let now $X$ be a smooth space. Assume that $\|e_i+e_j\|=1$ for some $1 \le i < j \le n$. Then 
$e^*_i(e_i+e_j)=e^*_j(e_i+e_j)=1$, hence the normalized support functional of the vector $e_i+e_j$ 
is not unique, which contradicts the smoothness of $X$. So $\|e_i+e_j\|>1$ and
similarly $\|e_i-e_j\|>1$. The rest of the proof proceeds as in the strictly convex case.

\end{proof}

\noindent \textbf{Remarks 2}\\
\noindent (1) It is clear that the assumption of strict convexity or smoothness
in Prop.3 can be replaced by $\|e_i \pm e_j\|>1$ for $1 \le i<j \le n$.\\
\noindent (2) An obvious example realizing Prop.3 is the set $\{\pm e_i:1 \le i \le n\}$ in
$\ell_p^n, \, 1<p<\infty, n \ge 2$.\\
\noindent (3) If $X$ is a Minkowski plane and $x,y \in S_X$
with $\|x-y\|=\|x+y\|=1$ then the points $\pm (x+y)$ and $\pm (x-y)$ lie
on $S_X$ and are the vertices of a parallelogram inscribed in the unit circle.
Since $\|\pm x\|=\|\pm y\|=1$, by Lemma 5, p.8 of \cite{SG} we get that all the segments 
joining neighbouring vertices must lie on the unit circle and hence the unit circle itself
coincides with the parallelogram with vertices $\pm (x+y),\pm (x-y)$. It follows
that $X$ is isometric to $\ell^2_\infty$. So if $X$ is not isometric to
$\ell^2_\infty$ and we have $x,y \in S_X$ with $\|x-y\|=1$, then necessarily $\|x+y\|>1$.
Similarly when $x,y \in S_X$ with $\|x-y\|=\|x+y\|=2$, it follows
that $X$ is isometric to $\ell^2_\infty$.

\begin{prop}
Let $(X,\| \cdot\|)$ be a Minkowski plane. Then $H'_\alpha(X)=4$.
\end{prop}

\beg
Let $\{(e_i,e^*_i):i=1,2\}$ be an Auerbach basis of $X$. Since the vectors $e_1$ and $e_2$
are mutually orthogonal, the unit ball $B_X$ is supported at $e_1$ by a line 
$L_1=\{z \in X:e^*_1(z)=1\}$ parallel to $e_2$ and also supported at $e_2$ by a line 
$L_2=\{z \in X:e^*_2(z)=1\}$ parallel to $e_1$. The lines $L_1,L_2$ are (non-parallel)
sides of the parallelogram with vertices $\{\pm(e_1-e_2),\pm(e_1+e_2)\}$ 
(and $\pm e_1, \pm e_2$ are the midpoints of these sides). Since $\{e_1,e_2\}$ is an
Auerbach basis, we get that $1 \le \|e_1 \pm e_2\|$ (and of course $\|e_1 \pm e_2\| \le 2$).
If all the vertices of this parallelogram are of norm 2, then $X$ is 
isometric to $\ell^2_\infty$ (see Remark 2(3)) and the result follows.

\begin{figure}
\centering

\includegraphics[width=0.75\textwidth]{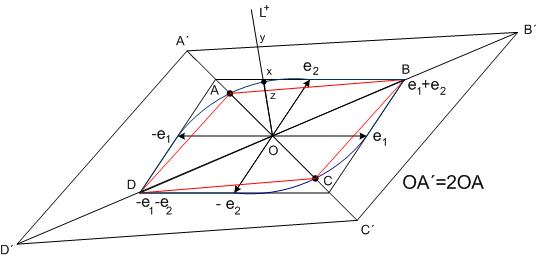}\\
\caption{}
\label{fig:1}

\end{figure}

Otherwise there is a pair of opposite vertices of norm less than 2 (say $\|\pm(e_1+e_2)\|<2$).
Now draw the diagonals of this parallelogram and take the points of intersection of these
diagonals with the unit sphere $S_X$ (see Fig.1). We get 4 points (the points $\pm \frac{e_1-e_2}{\|e_1-e_2\|},
\pm \frac{e_1+e_2}{\|e_1+e_2\|}$) which constitute a new parallelogram $ABCD$. This is a \bsa
set of 4 points lying on the unit sphere. Consider the norm one functionals $f$ and $g$
whose kernels are lines parallel to $AB$ and $BC$.

Let $A'B'C'D'$ be the homothetic copy of $ABCD$ of factor 2 with respect to the origin.
Since $A'B'C'D'$ does not touch the unit sphere, $f$ and $g$ give an evaluation
$>\frac{1}{2}$ and $<-\frac{1}{2}$ on opposite sides of $ABCD$. Indeed, let $x \in S_X$ such that
$f(x)=1$, then $x \notin$ Ker$f$, so the half-line $L^{+}=\{\lambda x:\lambda \ge 0\}$ either intersects
the line $A'B'$ or the line $D'C'$. Let for instance $L^{+}$ intersect $A'B'$ at $y=\lambda x$
and $AB$ at $z$. Then clearly $\lambda>1$ and $z=\frac{\lambda}{2} x$, thus $f(z)=\frac{\lambda}{2}
f(x)=\frac{\lambda}{2}>\frac{1}{2}$. Similarly we get that also $g$ has the desired property.
This way we obtain the constant $d>1$ in the definition of \bsa subset of $X$.
\end{proof}

\noindent \textbf{Remarks 3}\\
\noindent (1) It is easy to see that in any Minkowski plane $(X,\| \cdot \|)$,
there is an Auerbach basis $\{(e_i,e^*_i):i=1,2\}$ of $X$ such that $\|e_1-e_2\|>1$ and $\|e_1+e_2\|>1$.
Since any Auerbach basis of a 2-dimensional Banach space is 1-suppression unconditional, 
Prop.3 and Remark 2(1) imply that the set $A=\{\pm e_1,\pm e_2\}$ is \bsa with
$d>1$, so we get still another proof of Prop.4. It is not known to us if in any $n$-dimensional
Minkowski space ($n \ge 3$) there is an Auerbach basis satisfying $\|e_i\pm e_j\|>1$, $\forall i \neq j$.

In order to prove the 2-dimensional case, one can consider the inscribed centrally symmetric parallelogram of maximum area.
If we call two non-opposite vertices $e_1$ and $e_2$, and if $e_1+e_2$ is on the unit circle, then,
by area maximality, the line through $e_2$ parallel to $e_1$ supports the unit circle, so the segment
from $e_2$ to $e_1+e_2$ lies on the boundary of the unit ball and then the parallelogram with vertices
$\pm e_1$ and $\pm(e_1+e_2)$ also has maximum area. It follows that the segment from $e_1$ to $-e_2$
is also on the unit circle (because of a supporting line parallel to $e_1+e_2$) and the unit ball
is an affine regular hexagon. We may now pick any vertice $\alpha$ in the interior of the segment
from $e_1$ to $-e_2$; then the parallelogram with vertices $\pm \alpha$ and $\pm(e_1+e_2)$ (also of maximum area)
yields a basis with the desired property.

\noindent (2) If $X$ is an $n$-dimensional space with $n \ge 3$ then $H'_\alpha(X) \ge 4$, since any
Minkowski plane admits a \bsa set establishing this fact (see Prop.4).
Moreover one can obtain an equilateral set of 4 points yielding the same result. By a result 
of V.V. Makeev, in any 3-dimensional space there are 4 equidistant points which are also equidistant 
from their common barycenter (see \cite{MAK}). Assuming that the distance of a vertice
from the barycenter is 1 (and that the barycenter coincides with the origin) we obtain a 
$\lambda$-equilateral set of 4 points lying on the unit sphere with $\lambda>1$ (note that the extreme
points of the convex hull of a $\lambda$-equilateral set are exactly the points of the equilateral set and  
the common barycenter has distance $<\lambda$ from each of the extreme points). Taking into account 
Remarks 1(4) we have the result.

\section*{Antipodal Hadwiger number of $\ell_p$ spaces, $1<p \le \infty$ and the exponential growth of $H'_\alpha(X)$}

In this chapter we evaluate the antipodal and strict antipodal Hadwiger numbers of many (finite 
dimensional) $\ell_p$ spaces. Let $\alpha_n=\frac{\log n}{\log2}=\log_2 n,\, n \ge 2$.
This sequence is strictly increasing: $1=\frac{\log2}{\log2}=\alpha_2<\frac{\log3}{\log2}=\alpha_3<\frac{\log4}{\log2}=\alpha_4=2
< \dots<\alpha_n<\dots$.

\begin{prop}
Let $n \ge 2$ and $1< p \le \infty$, then the following hold:
\begin{enumerate}
\item When $p>\alpha_n$, then $H'_\alpha(\ell_p^n)=2^n$ (hence also $H_\alpha(\ell_p^n)=2^n$).
\item When $p \ge \alpha_n$, then $H_\alpha(\ell_p^n)=2^n$.
\end{enumerate} 
\end{prop}

\beg
Consider the set $S=\{-1,1\}^n$ and note that it is an antipodal subset of $\mathbb{R}^n$. 

Assume first that $1<p <\infty$. Then $\|x\|_p=n^{1/p}$, for every $x \in S$. Since 
$\lim_{p \rightarrow \infty} n^{1/p}=1$ and $1<n^{1/p}, \, n \ge 2$ and $p>1$, we get that there
is $p_0(n)$ such that \[p>p_0(n) \Rightarrow 1<n^{1/p}<2\;\;\;\;\;\;\;\; (1).\]
We will show that the least number $p_0(n)$ such that inequality (1) holds is $\alpha_n$. Indeed, 
for $n \ge 2$ and $p>1$ we have
\[n^{1/p}<2 \Leftrightarrow \log n^{1/p}<\log 2\Leftrightarrow \frac{1}{p} \log n<\log 2 \Leftrightarrow p> \alpha_n=
\frac{\log n}{\log 2}. \]

Set $S'=\frac{1}{n^{1/p}} \cdot S$, for $p \ge \alpha_n$ and $n \ge 2$ and observe that 
$S'$ is a \bsa subset of $S_{\ell_p^n}$
with constant $d=\frac{2}{n^{1/p}} \ge 1$. This is true because, given $x \neq y \in S'$, 
$x=(x_1,x_2,\dots,x_n),\, y=(y_1,y_2,\dots,y_n)$, there is $1 \le k \le n$ with $x_k \neq y_k$
such that the numbers $x_k,y_k$ have different signs and $|x_k|=|y_k|=\frac{1}{n^{1/p}}$.
Without loss of generality, let $x_k>0$ and $y_k<0$. Then
\[e^*_k(y)=-\frac{1}{n^{1/p}} \le e^*_k(z) \le e^*_k(x)=\frac{1}{n^{1/p}},\,\,\textrm{for} \,\, z \in S'.\]
It is clear from these inequalities that $H'_\alpha(\ell_p^n)=2^n$ for $p \in (\alpha_n, \infty)$ (since then
$d=\frac{2}{n^{1/p}} > 1$) and $H_\alpha(\ell_p^n)=2^n$ for $p=\alpha_n$ (since then $n^{1/p}=2$ and thus $d=1$).

In case when $p= \infty$, it is easily verified that the set $S$ itself is \bsa subset of $S_{\ell_\infty^n}$ with constant 
$d=2$ (actually $S$ is a 2-equilateral set, cf. Remark 1(4)), hence $H'_\alpha(\ell_\infty^n)=2^n$.

The proof of the Proposition is now complete.
\end{proof}

\noindent \textbf{Remarks 4}\\
\noindent (1) For $n=2$ we already have the stronger result $H'_\alpha(X)=4$ for any Minkowski plane, by Proposition 4.\\
\noindent (2) For $n=3$ we get that, when $p>\alpha_3=\frac{\log 3}{\log 2} \simeq 1.58$, then $H'_\alpha(\ell_p^3)=2^3=8$.
In particular $H'_\alpha(\ell_2^3)=2^3=8$.

The following Proposition settles the situation in case of a 3-dimensional space for the remaining $1<p \le \alpha_3 \simeq 1.58$:

\begin{prop}
When $1<p<\infty$, then $H'_\alpha(\ell_p^3)=2^3=8$.
\end{prop}

\beg
Let $2<\beta<2^p$, it then holds that $\frac{2}{\beta^{1/p}}>1$. Take $\alpha>0$
(and $\alpha<2$) such that $\frac{1}{\alpha}+\frac{1}{\beta}=1 \Leftrightarrow \alpha=\frac{\beta}{\beta-1}$.
We consider the points
\[x_1=\left(\frac{1}{\alpha^{1/p}},0,\frac{1}{\beta^{1/p}}\right) \,\,\,x_2=\left(-\frac{1}{\alpha^{1/p}},0,\frac{1}{\beta^{1/p}}\right)\]
\[x_3=\left(0,\frac{1}{\alpha^{1/p}},\frac{1}{\beta^{1/p}}\right) \,\,\,x_4=\left(0,-\frac{1}{\alpha^{1/p}},\frac{1}{\beta^{1/p}}\right).\]

Also set $x_5=-x_1,\,x_6=-x_2,\,x_7=-x_3$ and $x_8=-x_4$. 
Since $x_4-x_1=\left(-\frac{1}{\alpha^{1/p}},-\frac{1}{\alpha^{1/p}},0\right)=x_2-x_3$, the points $x_1,x_2,x_3,x_4$ are vertices of a
parallelogram (actually an orthogonal parallelogram). Thus, the points $\pm x_1,\pm x_2,\pm x_3,\pm x_4$ are vertices of an orthogonal
parallelepiped. Observe that $\|x_k\|_p=1$, $k=1,2,3,4$, so $\pm x_k \in S_{\ell_p^3}$, $k=1,2,3,4$ and also $\|x_k-x_l\|_p>1$,
$1 \le k<l \le 8$.

Obviously $S=\{x_k:k=1,2,\dots,8\}$ is an antipodal subset of $S_{\ell_p^3}$. We will show that it is \bsa with constant $d>1$.
Consider the functionals 
\[f_1=e_3^*,\,f_2=\frac{e_1^*+e_2^*}{2^{1/q}}\,\, \textrm{and}\,\,f_3=\frac{e_1^*-e_2^*}{2^{1/q}},\] 
where $q$ satisfies $\frac{1}{p}+\frac{1}{q}=1$ and $\{e_1,e_2,e_3\}$ is the usual basis of $\ell_p^3$, $1<p < \infty$. 
Note that $f_1,f_2,f_3 \in S_{\ell_q^3}$ (the unit sphere of the dual space $\ell_q=\ell_p^*$).
Geometrically, the kernels of $f_1,f_2,f_3$ are the planes $z=0,y=-x$ and $y=x$ of $\mathbb{R}^3$ respectively.

The following are easy to check:
\[f_1(x_1)=f_1(x_2)=f_1(x_3)=f_1(x_4)=\frac{1}{\beta^{1/p}} \,\, \textrm{and}\]
\[f_1(x_5)=f_1(x_6)=f_1(x_7)=f_1(x_8)=-\frac{1}{\beta^{1/p}}.\]

Therefore the points $\{x_1,x_2,x_3,x_4\}$ and $\{x_5,x_6,x_7,x_8\}$ are separated by planes parallel to $z=0$ and the difference
is $\frac{1}{\beta^{1/p}}-\left(-\frac{1}{\beta^{1/p}}\right)=\frac{2}{\beta^{1/p}}>1$. Also
\[f_2(x_1)=f_2(x_3)=f_2(x_6)=f_2(x_8)=\frac{1}{2^{1/q}} \cdot \frac{1}{\alpha^{1/p}}>
\frac{1}{2^{1/q}} \cdot \frac{1}{2^{1/p}}=\frac{1}{2^{1/p+1/q}}=\frac{1}{2} \]
\[\textrm{and} \,\, f_2(x_2)=f_2(x_4)=f_2(x_5)=f_2(x_7)=-\frac{1}{2^{1/q}} \cdot \frac{1}{\alpha^{1/p}}<-\frac{1}{2}.\]

So the separation now is achieved by planes parallel to $y=-x$ and the difference is $\frac{1}{2^{1/q}} \cdot \frac{1}{\alpha^{1/p}}-
\left(-\frac{1}{2^{1/q}} \cdot \frac{1}{\alpha^{1/p}}\right)=\frac{2}{2^{1/q} \cdot \alpha^{1/p}}>1$. We also have
\[f_3(x_1)=f_3(x_4)=f_3(x_6)=f_3(x_7)=\frac{1}{2^{1/q}} \cdot \frac{1}{\alpha^{1/p}}>\frac{1}{2} \,\, \textrm{and}\]
\[f_3(x_2)=f_3(x_3)=f_3(x_5)=f_3(x_8)=-\frac{1}{2^{1/q}} \cdot \frac{1}{\alpha^{1/p}}<-\frac{1}{2}.\]

The separation of the points is now achieved by planes parallel to $y=x$ and the corresponding difference is 
$\frac{1}{2^{1/q}} \cdot \frac{1}{\alpha^{1/p}}-\left(-\frac{1}{2^{1/q}} \cdot \frac{1}{\alpha^{1/p}}\right)=
\frac{2}{2^{1/q} \cdot \alpha^{1/p}}>1$.

From the above calculations, we conclude that $d=\min\{\frac{2}{\beta^{1/p}},\frac{2}{2^{1/q} \cdot \alpha^{1/p}}\}>1$
and thus $H'_\alpha(\ell_p^3)=2^3=8,\,\, \forall p>1$.

\end{proof}

Using the same method of proof, one can prove the following generalization:

\begin{prop}
When $n \ge 3$ and $1<p<\infty$, then $H'_\alpha(\ell_p^n) \ge 4n-4$.
\end{prop}

\beg
Assigning to the last coordinate the value $\pm \frac{1}{\beta^{1/p}}$ and placing
$\pm \frac{1}{\alpha^{1/p}}$ each time in one of the first $n-1$ coordinates
with 0 in every other place, one obtains the $4n-4$ required vectors.
\end{proof}

\noindent \textbf{Remark 5}\\
\noindent It is clear from Remark 4(1) that the above result also holds true for $n=2$. 
Since the Banach space $\ell_p^n,\,1<p<\infty,\, n \ge 2$ is smooth, we get from Prop. 2 that 
$H'_\alpha(\ell_p^n) \ge 2n$. Proposition 7 though provides us with a better
lower bound, $H'_\alpha(\ell_p^n) \ge 4n-4$.

\begin{prop}
Let $n \ge 4$  ($\alpha_n \ge 2$), then the following hold:
\begin{enumerate}
\item When $2 \le p \le \alpha_n$, then $4n-4 \le H'_\alpha(\ell_p^n)<2^n$.
\item When $2 \le p <\alpha_n$, then $4n-4 \le H'_\alpha(\ell_p^n) \le H_\alpha(\ell_p^n)<2^n$.
\item When $p=\alpha_n$, then $4n-4 \le H'_\alpha(\ell_p^n) < H_\alpha(\ell_p^n)=2^n$
\end{enumerate}
\end{prop}

\beg
To prove (1), given $p,n$ satisfying $2 \le p \le \alpha_n$ assume that the contrary holds,
i.e. $H'_\alpha(\ell_p^n)=2^n$. Then $B_{\ell_p^n}$ contains a \bsa subset $S$ with
constant $d>1$ and cardinality $|S|=2^n$. By Theorem 3 of \cite{MV} stated in the 
Introduction, there is an equivalent norm $||| \cdot |||$ on $\mathbb{R}^n$ which admits
an equilateral set of cardinality $2^n$, hence $(\mathbb{R}^n,|||\cdot|||)$ is isometric to
$\ell_\infty^n$. Moreover the same Theorem yields for the Banach-Mazur distance of the norms
$\|\cdot\|_p$ and $|||\cdot|||$ that
\[d(\ell_p^n,\ell_\infty^n) \le \frac{2}{d}<2,\,\, \textrm{since}\,\,d>1.\]
But for $p \ge 2$ we know that $d(\ell_p^n,\ell_\infty^n)=n^{1/p}$, see \cite{FHHMZ}.
Hence \[n^{1/p}=d(\ell_p^n,\ell_\infty^n) \le \frac{2}{d}<2 \Rightarrow p>\frac{\log n}{\log 2}=\alpha_n \]
which contradicts our assumption. Taking into account Proposition 7, the proof of (1) is complete.

The proof of (2) is similar. Concerning (3), the inequalities follow from (1) and the equality
follows from Proposition 5(2).
\end{proof}

The following Theorem summarizes all the previous results about $\ell_p^n$, $1 <p \le +\infty$ spaces:

\begin{theo}
Let $1 <p \le +\infty$, then the following hold:
\begin{enumerate}
\item When $n=2$ or $3$ , then $H'_\alpha(\ell_p^n)=4n-4=2^n=H_\alpha(\ell_p^n)$.
\item When $n \ge 4$, then we have:
\begin{enumerate}
\item $H'_\alpha(\ell_p^n) \ge 4n-4$
\item when $2 \le p<\alpha_n$, then $4n-4 \le H'_\alpha(\ell_p^n) \le H_\alpha(\ell_p^n)< 2^n$
\item when $p=\alpha_n$, then $4n-4 \le H'_\alpha(\ell_p^n)<2^n=H_\alpha(\ell_p^n)$, in particular
for $n=4$ we get that $p=\alpha_4=2$ and $12 \le H'_\alpha(\ell_2^4)<16=H_\alpha(\ell_2^4)$ and
\item when $p>\alpha_n$, then $H'_\alpha(\ell_p^n)=2^n=H_\alpha(\ell_p^n)$.
\end{enumerate}
\end{enumerate}
\end{theo}

We now introduce the number $e_c(X)$ for a finite-dimensional Banach space $X$ and show that it is 
bounded below by an unbounded function of the dimension of $X$. Since $H'_\alpha(X) \ge e_c(X)$,
we have that $H'_\alpha(X)$ also has this property.

Recall that $e(X)$ (=the equilateral number of $X$)
denotes the largest size of an equilateral subset of $X$. We define analogously the number $e_c(X)$
to be the largest size of a $\lambda$-equilateral subset of $S_X$, where $1<\lambda \le 2$. 
It is clear that
\[e_c(X) \le e(X)\]
and also by Remark 1(4) that
\[H'_\alpha(X) \ge e_c(X).\]
Now, one can prove by similar arguments (we omit the details) the following variant of a significant
Theorem of Brass and Dekster (see Theorem 8 of \cite{S}):

\begin{theo}
Let $X$ be an $n$-dimensional Banach space ($n \ge 2$) with Banach-Mazur distance $d(X,\ell_2^n) 
\le 1+ \frac{1}{3(n+1)}$. Then any $\lambda$-equilateral set in $S_X$, where $\lambda \in \left(
1,1+ \frac{1}{3(n+1)} \right)$, of at most $n-1$ points can be extended to a $\lambda$-equilateral set 
in $S_X$ of $n$ points.
\end{theo}

We can prove using Dvoretzky's Theorem and Theorem 2, in the same way as Theorem 7 of \cite{S}
is proved, that if dim$X=n$ then $e_c(X) \ge \kappa(logn)^\frac{1}{3}$ for some constant $\kappa>0$ and $n$
sufficiently large.

So from the above observations we get the following:

\begin{theo}
Let $X$ be an $n$-dimensional Banach space. Then \[H'_\alpha(X) \ge e_c(X) \ge \kappa(logn)^\frac{1}{3}.\]
\end{theo}

\noindent \underline{Note.} Using the techniques of Swanepoel and Villa in \cite{SV} (see also Th.4.3 of \cite{K})
one can show that $H'_\alpha(X) \ge e_c(X) \ge e^{\kappa_1 \sqrt{logn}}$, for some constant $\kappa_1>0$. We also 
note that for $n=2$ or $3$ the inequality $e_c(X) \ge n+1$ holds; for the case $n=2$ we refer the reader 
to \cite{KO}, Prop. 1.2 and for $n=3$ to Remarks 3(2). It is not difficult to check that $e_c(\ell_p^n) \ge n+1$ 
for $1 < p <+\infty$, $n \ge 1$ and in particular that $e_c(\ell_2^n)=n+1$. Also (obviously) 
$e_c(\ell_1^n) \ge 2n$ and $e_c(\ell_{\infty}^n)=2^n$.

In the sequel we will prove that $H'_\alpha(X)$ actually increases exponentially in the dimension of $X$. We
need the following Lemma (see \cite{S}, Lemma 2), which according to \cite{S} is a special
case of the Johnson-Lindenstrauss flattening Lemma and also results from the Gilbert-Varshamov
lower bound for binary codes, see also \cite{RS} Th.13.

\begin{lemm}
For each $\delta>0$ there exist $\varepsilon=\varepsilon(\delta)$ and $n_0=n_0(\delta) \ge 1$
such that for all $n \ge n_0$ there exist $m>(1+\varepsilon)^n$ vectors 
$w_1,\cdots,w_m \in \mathbb{R}^n$ ($w_i=\frac{x_i}{\sqrt{n}},\, x_i \in \{-1,1\}^n$) satisfying

\begin{displaymath}
\left\{ \begin{array}{ll}
<w_i,w_i>=1, & i=1,2,\cdots,m\\
|<w_i,w_j>|<\delta, & \textrm{for all distinct } i,j \in \{1,2,\cdots,m\}\\
\end{array} \right.
\end{displaymath}
We may take $\varepsilon=\frac{\delta^2}{2}$ and $n_0 \ge \frac{120log2}{25 \delta^4-\delta^6}$.
\end{lemm}

We first prove the result for the spaces $\ell_2^n$. The essential part of this proof
is contained in the proof of Th.2 in \cite{S}.
\begin{prop}
Let $X=\ell_2^n$, then the number $H'_\alpha(X)$ increases exponentially in the dimension of $X$.
\end{prop}

\beg Let $0<\delta \le \frac{1}{3}$, then by Lemma 1 
there exist $\varepsilon=\varepsilon(\delta)>0$ and $n_0=n_0(\delta) \ge 1$
such that for all $n \ge n_0$ there exist $m>(1+\varepsilon)^n$ vectors 
$w_1,\cdots,w_m \in \mathbb{R}^n$ satisfying the two conclusions. We will prove
that the set $\{w_1,\cdots,w_m\}$ is \bsa subset of $S_{\ell_2^n}$ with
constant $d>\sqrt{2-2\delta}\,( \ge \sqrt{\frac{4}{3}}>1)$.

By the first conclusion of Lemma 1 we have $||w_i||_2=1$ for $i=1,\cdots,m$. For any
$i \neq j \in \{1,2,\cdots,m\}$ we define the linear functional 
\[f(x)=<x,w_i-w_j>,\, x \in \mathbb{R}^n\] 
with $||f||_2=||w_i-w_j||_2$.

Let $k \in \{1,2,\cdots,m\}$ with $k \neq i,j$, we will prove that
\[f(w_j) \le f(w_k) \le f(w_i) \,\,\,(1).  \]
It suffices to prove that\\ 
$\max\{|f(w_i-w_k)|,|f(w_j-w_k)|:k=1,2,\cdots,m \text{ and }
k \neq i,j\} \le f(w_i-w_j)=||w_i-w_j||_2^2$.
Setting $\alpha=<w_i,w_j>$
(and taking into account Lemma 1 and the fact that $0<\delta \le \frac{1}{3}$) we have
$f(w_i-w_k)=<w_i-w_k,w_i-w_j>=
<w_i,w_i>-<w_i,w_j>-<w_k,w_i>+<w_k,w_j>\,<1-\alpha+\frac{2}{3}=\frac{5}{3}-\alpha$.

Also $f(w_i-w_j)=<w_i-w_j,w_i-w_j>=
<w_i,w_i>+<w_j,w_j>$ $-2<w_i,w_j>=2-2\alpha$. Since $|\alpha|=|<w_i,w_j>|<\delta
\le \frac{1}{3}$, it holds that $\frac{5}{3}-\alpha<2-2\alpha$.
Similarly $-f(w_i-w_k)=\alpha-\frac{5}{3}<2-2\alpha$, hence $|f(w_i-w_k)|<f(w_i-w_j)$.
The other inequality required can be proved analogously and so (1) holds true.

Now setting $g=\frac{f}{||f||_2}$ we have $||g||_2=1$ and\\ 
$g(w_i-w_j)=\frac{1}{||f||_2}
f(w_i-w_j)=\frac{1}{||w_i-w_j||_2}||w_i-w_j||_2^2=||w_i-w_j||_2=\sqrt{2-2<w_i,w_j>}>
\sqrt{2-2\delta} \ge \sqrt{\frac{4}{3}}>1$ and the proof of the 
Proposition is complete. 
\end{proof}

\noindent \textbf{Remark 6}\\
The proof of Proposition 9 implies that, choosing $0<\delta \le \frac{1}{3}$ small enough, 
the distance $||w_i-w_j||_2=\sqrt{2-2<w_i,w_j>}>
\sqrt{2-2\delta}$ (for $i \neq j$) becomes arbitrarily close to $\sqrt{2}$. Also
note that the set $\{w_1,\cdots,w_m\}$ is strictly antipodal, since the inequalities (1)
of Prop.9 are strict.

Using the previous Proposition and an important result of Milman (the Quotient of Subspace Theorem, see
\cite{M}, also \cite{GIM})
the way it was used by Bourgain in \cite{FL}, Th. 4.3, we are now in a position to prove the 
general result:

\begin{theo}
Let $X$ be a finite-dimensional Banach space, then the number $H'_\alpha(X)$ 
increases exponentially in the dimension of $X$.
\end{theo}

\beg From Milman's result there is a function $\psi:(0,+\infty) \rightarrow (0,+\infty)$,
such that for all $\delta>0$ and for every Banach space $X$ with dim$X=n$, there exist
subspaces $Z \subseteq Y \subseteq X$ satisfying dim$Z \ge \psi(\delta)n$ and there is  
an ellipsoid $L \subseteq Z$, with
\[L \subseteq \pi(B_Y) \subseteq (1+\delta)L \;\;\;(1) \]
where $\pi:Y \rightarrow Z$ is the orthogonal projection of $Y$ onto $Z$. (There exists a proper
euclidean structure on $X$ used in the proof of Milman's Theorem. We consider the orthogonal 
projection $P_Z:Y \rightarrow Z$ with respect to this structure, which algebraically coincides
with the quotient map $\pi:Y \rightarrow Y/Z^{\perp}$; hence $\pi(B_Y)$ induces the quotient
norm $||z||_{qt}=\inf\{||z-z_0||:z_0 \in Z^{\perp}\}$ on $Z$ which is isomorphic with $Y/Z^{\perp}$,
see \cite{GIM} Ths. 5.2.1, 5.3.1 and \cite{VER}, Prop. 2.3).

Assume that $\delta \in \left(0, \frac{1}{3} \right)$. Then by the previous Proposition there is 
$\varepsilon=\varepsilon(\delta)>0$, such that for $n \in \mathbb{N}$ large enough there 
exist $m \ge (1+\varepsilon)^{\psi(\delta)n}$ vectors $z_1,z_2,\cdots,z_m$ on the boundary
of the euclidean unit ball $L$ in $Z$ with the set $S=\{z_1,z_2,\cdots,z_m\}$ being a \bsa
set with constant $d>\sqrt{2-2\delta}$. This means that for $i \neq j \in \{1,2,\cdots,m\}$
there is a linear functional $f \in Z^*$ with $||f||_{L}=1$ and such that
\[f(z_i)-f(z_j) \ge d \text{ and } f(z_j) \le f(z_k) \le f(z_i), \text{ for } 
1 \le k \le m. \;\;\;(2)\]
Now (1) is equivalent to
\[||z||_{qt} \le ||z||_L \le (1+\delta)||z||_{qt},\, z \in Z \;\;\;(3)\]
which implies for the dual norms that
\[||f||_L \le||f||_{qt} \le (1+\delta)||f||_L,\, f \in Z^*.\;\;\;(4)\]
It is clear that for $i \le m$ there is $y_i \in B_Y$ such that
$\pi(y_i)=z_i$. We note that we can select $y_i \in S_Y$ so that
$\pi(y_i)=z_i$ holds (if $z_i \in S_Y$ then set $y_i=z_i$,
otherwise for any $z_0 \in Z^{\perp} \setminus \{0\}$, there is $\lambda>0$ so that 
$y_i=z_i+\lambda z_0 \in S_Y$ and $\pi(y_i)=z_i$).

Next, we will prove that the set $S^\prime=\{y_1,y_2,\cdots,y_m\}$ is a \bsa subset of $S_Y$ (hence also of $S_X$)
with constant $\ge \frac{\sqrt{2-2\delta}}{1+\delta}$. Let $1 \le i \neq j \le m$. Take an $f \in Z^*$
with $||f||_L=1$ such that the inequalities in (2) hold and set $g=f \circ \pi$. Then (using (4) and $||\pi||=1$)
we have 
\[||g|| \le ||\pi|| \cdot ||f||_{qt} = ||f||_{qt} \le (1+\delta) ||f||_L \le 1+\delta \;\;\;(5).\]

Observe that, since $g(y_k)=f(\pi(y_k))=f(z_k) ,\,1 \le k \le m$, from (2) we have 
\[g(y_j) \le g(y_k) \le g(y_i), \, \text{ for } 1 \le k \le m.  \]
Also $g(y_i)-g(y_j)=f(\pi(y_i))-f(\pi(y_j))=f(z_i)-f(z_j) \ge d$.
Relation (5) implies that
\[\frac{g}{||g||}(y_i)-\frac{g}{||g||}(y_j)=\frac{1}{||g||}(f(z_i)-f(z_j)) \ge \frac{d}{1+\delta}>\frac{\sqrt{2-2\delta}}{1+\delta}. \]
For $0<\delta<\sqrt{5}-2$ we have that $\frac{\sqrt{2-2\delta}}{1+\delta}>1$
(and $\frac{\sqrt{2-2\delta}}{1+\delta} \rightarrow \sqrt{2}$ when $\delta \rightarrow 0,\, \delta>0$).
The proof is now complete.

\end{proof}

We also investigated the spaces $\ell_1^3$ and the Petty space on $\mathbb{R}^3$ (see \cite{S}, \cite{SWAN})
with respect to their (strict) antipodal Hadwiger numbers. In both spaces, the strict antipodal 
Hadwiger numbers are as big as possible:

\begin{prop}

\begin{enumerate}
\item  $H_\alpha(\ell_1^3)=H'_\alpha(\ell_1^3)=2^3=8$.
\item If $(X,\| \cdot \|)$ is the Petty space on $\mathbb{R}^3$, where $\|(x,y,z)\|=
\sqrt{x^2+y^2}+|z|$, then $H_\alpha(X)=H'_\alpha(X)=2^3=8$.
\end{enumerate}

\end{prop}

\beg 
To prove (1), set $x_1=\left(1,1,-\frac{1}{3}\right),\, x_2=\left(1,-\frac{1}{3},1\right),\,x_3=\left(-\frac{1}{3},1,1 \right)$ and 
$O=$conv $\{\pm x_1, \pm x_2,\pm x_3\}$. The Minkowski functional of $O$ defines a norm and the corresponding space is isometric to $\ell_1^3$
through the linear isometry designated by $T(e_i)=x_i,\, i=1,2,3$ ($\{e_1,e_2,e_3\}$ is the usual basis of $\ell_1^3$). Let also $C_3=B_{\ell_\infty^3}=[-1,1]^3$. In \cite{X} Fei Xue
observed that the octahedron $O$ and the cube $C_3$ satisfy
\[  \frac{5}{9} C_3 \subseteq O \subseteq C_3\]
and obtained for the Banach-Mazur distance of the spaces $\ell_1^3$ and $\ell_\infty^3$ the upper bound
$d(\ell_1^3,\ell_\infty^3) \le \frac{9}{5}<2$.

The set $S=\frac{5}{9} \{-1,1\}^3$ is the set of vertices of a parallelepiped and all of its points belong to
the boundary of the octahedron $O$ (hence are vectors of $\ell_1$-norm 1). One can readily check that
\[ \frac{5}{9} (1,1,1)=\frac{1}{3}(x_1+x_2+x_3),\; \;\; \frac{5}{9} (1,1,-1)=\frac{2}{3} x_1+\frac{1}{6} (-x_2)+\frac{1}{6} (-x_3) \]

\[ \frac{5}{9} (-1,1,1)=\frac{2}{3} x_3+ \frac{1}{6} (-x_1)+ \frac{1}{6} (-x_2),\;\;\; 
\frac{5}{9} (-1,1,-1)=\frac{2}{3} (-x_2)+\frac{1}{6} x_3+\frac{1}{6} x_1 \]
and the other points are the symmetric of these 4. We will show that $S$ is a \bsa subset of $S_{\ell_1^3}$ with constant 
$d=\frac{10}{9}>1$.

The functionals separating opposite faces of the parallelepiped are the $e_i^*,\, i=1,2,3$ (since $O \subseteq C_3$, we have that
$|e_i^*(x,y,z)| \le 1$ for any $(x,y,z) \in O$, hence $e_i^* \in B_{(\ell_1^3)^*}$).
For instance 
\[e_1^*\left(\frac{5}{9},\frac{5}{9},\frac{5}{9}\right)=e_1^*\left(\frac{5}{9},-\frac{5}{9},\frac{5}{9}\right)=
e_1^*\left(\frac{5}{9},\frac{5}{9},-\frac{5}{9}\right)=e_1^*\left(\frac{5}{9},-\frac{5}{9},-\frac{5}{9}\right)=\frac{5}{9}\;\; \textrm{while}\]

\[e_1^*\left(-\frac{5}{9},\frac{5}{9},\frac{5}{9}\right)=e_1^*\left(-\frac{5}{9},-\frac{5}{9},\frac{5}{9}\right)=
e_1^*\left(-\frac{5}{9},\frac{5}{9},-\frac{5}{9}\right)=e_1^*\left(-\frac{5}{9},-\frac{5}{9},-\frac{5}{9}\right)=-\frac{5}{9}\]

and the evaluations for the other faces are similar (see also the proof of Proposition 6), which implies that $d=\frac{5}{9}-(-\frac{5}{9})=\frac{10}{9}>1$.

Now for the proof of (2) one may consider the points $A(-0.18,0,0.82)$, $B(0.82,0,-0.18)$, $C(0.32,0.6,0.32)$ and
$D(0.32,-0.6,0.32)$. These points lie on the unit sphere of the Petty space and form a parallelogram, as $\overrightarrow{AC}=
\overrightarrow{DB}=(0.5,0.6,-0.5)$. So these points together with the symmetric points $A',B',C',D'$ with respect to the origin
form a parallelepiped with all vertices on the unit sphere. To find the three functionals separating opposite faces, we first
calculate the equations of three planes, each one defined by three vertices of the parallelepiped. We have:
\[ADB:x+z=0.64\]
\[ADC':0.6x+y-0.6z=-0.6 \;\; \textrm{and}\]
\[C'DB:0.6x-y-0.6z=0.6.\]
Taking into account the dual norm $\|(x,y,z)\|_*=\max\{\sqrt{x^2+y^2},|z|\}$ and normalizing, each of these planes 
yields a functional which, of course, evaluates all points on a face the same way. So set $f_1=(1,0,1),\,f_2=\frac{1}{\sqrt{1.36}}(0.6,1,-0.6)$ and
$f_3=\frac{1}{\sqrt{1.36}}(0.6,-1,-0.6)$ which are all of norm 1. The separation of the faces goes as follows:
\[f_1(A)=f_1(D)=f_1(B)=f_1(C)=0.64\]
\[f_1(A')=f_1(D')=f_1(B')=f_1(C')=-0.64\]
with $f_1(A)-f_1(A')=1.28>1$.
\[f_2(A)=f_2(D)=f_2(B')=f_2(C')=-\frac{0.6}{\sqrt{1.36}}\]
\[f_2(A')=f_2(D')=f_2(B)=f_2(C)=\frac{0.6}{\sqrt{1.36}}\]
also 
\[f_3(A')=f_3(D)=f_3(B)=f_3(C')=\frac{0.6}{\sqrt{1.36}}\]
\[f_3(A)=f_3(D')=f_3(B')=f_3(C)=-\frac{0.6}{\sqrt{1.36}}\]
with $f_2(A')-f_2(A)=f_3(A')-f_3(A)=\frac{1.2}{\sqrt{1.36}}\simeq 1.02899>1$. 
Hence $d \simeq 1.02899>1$ and the conclusion follows.

\end{proof}

\noindent \textbf{Remarks 7}\\
\noindent (1) Using the Theorem stated in the Introduction (Th. 3 of \cite{MV}) we obtain an upper bound
for the Banach-Mazur distance of the Petty space $X$ from $\ell_\infty^3$. We apply the Theorem as
in the proof of Proposition 8 and conclude that $d(X,\ell_\infty^3)<2$.
One can find a better upper bound by direct calculation. Since $d(X,\ell_\infty^3)=d(X^*,\ell_1^3)$,
it suffices to evaluate the distance between the dual spaces, which is easier. The ball of $X^*$ is the
right cylinder $B_{X^*}=\{(x,y,z):\sqrt{x^2+y^2} \le 1\;\;\textrm{and}\;\;|z| \le 1\}$.
Consider the points $A(0,-1,1),\,B(0.8,0.6,1)$ and $C(-0.8,0.6,1)$ of $S_{X^*}$ along with the symmetric points
$A',B',C'$ with respect to the origin. Let $O=$ conv $\{\pm \overrightarrow{OA},\pm \overrightarrow{OB},\pm \overrightarrow{OC}\}$
(an octahedron producing a space linearly isometric to $\ell_1^3$). We calculate the largest $\alpha>\frac{1}{2}$
such that
\[\alpha \cdot B_{X^*} \subseteq O \subseteq B_{X^*}.\]
The right-hand inclusion is obvious. To find the optimal value of $\alpha$ one has to check which homothetic
copy of $B_{X^*}$ touches some of the faces of the octahedron in a single point (it suffices to check the upper 4 faces
due to symmetry). The upper 4 faces define the planes
\[ABC:z=1\]
\[CAB':10x+5y-3z+8=0\]
\[BAC':10x-5y+3z-8=0 \;\; \textrm{and}\]
\[AB'C':5y+z=-4.\]
Note that if the upper half of the octahedron touches a homothetic copy $ \alpha \cdot B_{X^*}$ of the cylinder, 
then the common point must lie on the circle $x^2+y^2=\alpha^2,\,z=\alpha$. Solving
the systems of equations of each of the above 4 planes together with $x^2+y^2=\alpha^2$ and $z=\alpha$ 
we find that for $\alpha\simeq 0.56416$ there is at most one common point of a face with the corresponding 
homothetic copy of $B_{X^*}$, hence $0.56416  \cdot B_{X^*} \subseteq O \subseteq B_{X^*}$ and thus 
$d(X,\ell_\infty^3)=d(X^*,\ell_1^3) \le \frac{1}{0.56416} \simeq 1.77254$.\\
\noindent (2) Concerning the original Hadwiger number of the Petty space $X$, we obtain a lower bound of 14
by using the 1-separated set of 14 points on the unit sphere of $X$:
\[\alpha_1=e_1=(1,0,0), \,\alpha_2=\left( \frac{1}{2},\frac{2}{3},\frac{1}{6} \right), \, 
\alpha_3=\left( 0,\frac{2}{3},-\frac{1}{3} \right),\,\alpha_4=\left( -\frac{1}{2},\frac{2}{3},\frac{1}{6} \right),\]
\[\alpha_5=-e_1,\,\alpha_6=\left(-\frac{1}{2},-\frac{2}{3},\frac{1}{6} \right),\,\alpha_7=\left( 0,-\frac{2}{3},-\frac{1}{3} \right),\, 
\alpha_8=\left( \frac{1}{2},-\frac{2}{3},\frac{1}{6} \right),\] 
\[b_1=e_3=(0,0,1), \, b_2=\left( \frac{1}{2},0,\frac{1}{2}\right),\, b_3=\left( -\frac{1}{2},0,\frac{1}{2}\right),\, b_4=-e_3,\]
\[b_5=\left( \frac{1}{2},0,-\frac{1}{2}\right),\, b_6=\left( -\frac{1}{2},0,-\frac{1}{2}\right).  \]
Observe that the points $\alpha_1,\alpha_2$ are the symmetric of $\alpha_5,\alpha_4$ with respect to the YZ-plane,
the points $\alpha_8,\alpha_7,\alpha_6$ are the symmetric of $\alpha_2,\alpha_3,\alpha_4$ with respect to the XZ-plane
and the points $b_1,b_2,b_3$ are the symmetric of $b_4,b_5,b_6$ with respect to the XY-plane.

Also the set of points
\[c_1=\left( \frac{2}{3},0,-\frac{1}{3}\right) \, c_2=\left( \frac{\sqrt{2}}{2},\frac{\sqrt{2}}{2},0 \right), \, 
c_3=\left( 0,\frac{2}{3},\frac{1}{3} \right),\]
\[c_4=\left( -\frac{\sqrt{2}}{2},\frac{\sqrt{2}}{2},0 \right),\,c_5=\left( -\frac{2}{3},0,-\frac{1}{3}\right), \,  c_6=\left( -\frac{\sqrt{2}}{2},-\frac{\sqrt{2}}{2},0 \right), \] 
\[c_7=\left( 0,-\frac{2}{3},\frac{1}{3} \right),\,c_8=\left( \frac{\sqrt{2}}{2},-\frac{\sqrt{2}}{2},0 \right),\,
b_1,\, b_2, \, b_3, \, b_4,\] 
\[b'_5=\left(0,\frac{1}{2},-\frac{1}{2}\right), \, b'_6=\left(0,-\frac{1}{2},-\frac{1}{2}\right),\] 
shows that $H(X) \ge 14$. Both of the above pointsets are maximal 1-separated subsets of $S_X$.
Furthermore the subset $\{c_1,c_2, \cdots,c_8\} \cup \{b_1,b_4\}$ of the second pointset yields the lower bound of 10
for the strict Hadwiger number of $X$, i.e. $H'(X) \ge 10$. It seems likely that $H(X)=14$, but the details
of such a proof are not yet clear.

\noindent (3) About the Hadwiger number of the euclidean spaces the exact values are known in case when $n=2,3,4,8$ and $24$
(see \cite{B}, also \cite{SWAN}).
When $2 \le n \le 6$, the best-known lower bounds are larger than $2^n$ and for $n=7$ we have $126 \le H(\ell_2^7) \le 134$.
Unlike that, when $8 \le n \le 24$, the best-known upper bounds are smaller than $2^n$.

We conclude with some open problems and questions:
\begin{enumerate}
\item Find better upper and lower bounds for the numbers:
\begin{enumerate}
\item $H'_\alpha(\ell_p^n)$ and $H_\alpha(\ell_p^n)$, for $n \ge 4$ and $1 \le p<2$,
\item $H'(\ell_{\alpha_n}^n)$, for $n \ge 4$ and
\item $H'_\alpha(\ell_2^n)$ and $H_\alpha(\ell_2^n)$, for $5 \le n \le 24$
\end{enumerate}
Of particular interest is the case $n=4$ (then $\alpha_4=2$)
and $n=8,\, 24$, where the exact values $H(\ell_2^n)$ are known.
\item Is there a 3-dimensional Banach space $X$, with $H_\alpha(X)<8$? Such a space would
be a candidate to have Banach-Mazur distance $d(X,\ell_\infty^3) \ge 2$ (see the proof of Prop.8).
\item Let $X$ be an $n$-dimensional Banach space with $n \ge 4$. Does the inequality 
$e_c(X) \ge n+1$ (or at least $H'_\alpha(X) \ge n+1$) hold?
(See Theorem 3 and the Note following it).
\end{enumerate}

\noindent S.K.Mercourakis, G.Vassiliadis\\
University of Athens\\
Department of Mathematics\\
15784 Athens, Greece\\
e-mail: smercour@math.uoa.gr

\hspace{0.3cm} georgevassil@hotmail.com
\end{document}